\documentclass[upref]{amsart}
\usepackage{amssymb}
\usepackage{times}

\numberwithin{equation}{section}

\newcommand{\lab}[1]{\label{#1}}

\newcommand{\ndash}{\nobreakdash-}
\newcommand{\Ndash}{\nobreakdash--}

\newtheorem{Thm}{Theorem}

\theoremstyle{remark}
\newtheorem{Rem}{Remark}[section]
\newtheorem*{Rem*}{Remark}
\newtheorem*{Ack}{Acknowledgment}

\theoremstyle{definition}

\newtheorem{Ass}{Assumption}

\newtheorem*{Conv*}{Convention}

\newcommand{\Lie}[2]{{\left[{\,{#1}\,,\,{#2}\,}\right]}}

\newcommand{\braket}[2]{\left\langle{\,{#1}\,,\,{#2}\,}\right\rangle}
\newcommand{\dbraket}[2]{\left\langle\!\left\langle
{\,{#1}\ ;\,{#2}\,}\right\rangle\!\right\rangle}

\newcommand{\BV}[2]{\left({\,{#1}\,,\,{#2}\,}\right)}
\newcommand{\dLie}[2]{{\left[\!\left[{\,{#1}\ ;\,{#2}\,}\right]\!\right]}}

\newcommand{\bbR}{{\mathbb{R}}}

\newcommand{\bbZ}{{\mathbb{Z}}}

\newcommand{\de}{\partial}

\newcommand{\SBV}{{S_{\mathrm{BV}}}}
\newcommand{\ii}{{\mathrm{i}}}

\newcommand{\dleft}[1]{\,\frac{\overset\leftarrow\de}
{\de{#1}}\,}

\newcommand{\dright}[1]{\,{\frac{\overset\rightarrow\de}
{\de{#1}}}\,}

\newcommand{\calA}{\mathcal{A}}

\newcommand{\calH}{\mathcal{H}}

\newcommand{\calG}{\mathcal{G}}

\newcommand{\bmu}{\boldsymbol{\mu}}
\newcommand{\blambda}{\boldsymbol{\lambda}}

\newcommand{\sfcalH}{\boldsymbol{\mathcal{H}}}

\newcommand{\sfA}{{\mathsf{A}}}
\newcommand{\sfB}{{\mathsf{B}}}
\newcommand{\sfF}{{\mathsf{F}}}
\newcommand{\sfa}{{\mathsf{a}}}
\newcommand{\sfh}{{\mathsf{h}}}

\newcommand{\sfO}{{\mathsf{O}}}
\newcommand{\sfS}{{\mathsf{S}}}

\newcommand{\frg}{{\mathfrak{g}}}
\newcommand{\frgl}{{\mathfrak{gl}}}

\newcommand{\dod}[1]{{\frac\dd{\dd{#1}}}}

\DeclareMathOperator{\Tr}{Tr}
\newcommand{\dd}{\mathrm{d}}

\DeclareMathOperator{\ad}{ad}
\DeclareMathOperator{\Ad}{Ad}
\DeclareMathOperator{\hol}{Hol}
\DeclareMathOperator{\sfhol}{{\mathsf{Hol}}}
\newcommand{\LL}{\mathrm{L}}
\DeclareMathOperator{\Imb}{Imb}
\newcommand{\dBRST}{{\delta_{\mathrm{BRST}}}}
\newcommand{\dBV}{{\delta_{\mathrm{BV}}}}
\newcommand{\DBV}{{\Delta_{\mathrm{BV}}}}
\newcommand{\OBV}{{\Omega_{\mathrm{BV}}}}
\DeclareMathOperator{\gh}{gh}
\newcommand{\sfdelta}{{\boldsymbol{\delta}}}
\newcommand{\sfDelta}{{\boldsymbol{\Delta}}}
\newcommand{\sfOmega}{{\boldsymbol{\Omega}}}

\newcommand{\ImbfM}{\Imb_f(S^1,M)}

\newcommand\qq{\rm}
\newcommand\cmp[1]{{\qq Commun.\ Math.\ Phys.\ \bf #1}}
\newcommand\jmp[1]{{\qq J.\ Math.\ Phys.\ \bf #1}}
\newcommand\pl[1]{{\qq Phys.\ Lett.\ \bf #1}}
\newcommand\np[1]{{\qq Nucl.\ Phys.\ \bf #1}}

\newcommand\lmp[1]{{\qq Lett.\ Math.\ Phys.\ \bf #1}}

\newcommand\ijmp[1]{{\qq Int.\ J. Mod.\ Phys.\ \bf #1}}

\newcommand\prept[1]{{\qq Phys.\ Rept.\ \bf #1}}

\newcommand\anp[1]{{\qq Ann.\ Phys.\ \bf #1}}
\newcommand\anm[1]{{\qq Ann.\ Math.\ \bf #1}}

\newcommand\phs[1]{{\qq Physica Scripta \bf #1}}

\begin{document}
\title[Loop observables for $BF$ theories in any dimension\dots]
{Loop observables for $BF$ theories in any dimension
and the cohomology of knots}

\author[A.~S.~Cattaneo]{Alberto~S.~Cattaneo}
\address{Mathematisches Institut --- Universit\"at Z\"urich--Irchel ---  
Winterthurerstrasse 190 --- CH-8057 Z\"urich --- Switzerland}  
\email{asc@math.unizh.ch}
\thanks{A.~S.~C. acknowledges partial support of SNF Grant No.~2100-055536.98/1}

\author[P.~Cotta-Ramusino]{Paolo~Cotta-Ramusino}  
\address{Dipartimento di Matematica --- Universit\`a degli Studi di 
Milano ---  via Saldini, 50 --- I-20133 Milano --- Italy; \&  
INFN, Sezione di Milano --- via Celoria, 16 --- I-20133 Milano --- Italy}  
\email{cotta@mi.infn.it}  

\author[C.~A.~Rossi]{Carlo~A.~Rossi}
\address{Mathematisches Institut --- Universit\"at Z\"urich--Irchel ---  
Winterthurerstrasse 190 --- CH-8057 Z\"urich --- Switzerland}  
\email{crossi@math.ethz.ch}


\date{}

\begin{abstract}
A generalization of Wilson loop observables for $BF$ theories in any dimension
is introduced in the Batalin--Vilkovisky framework. The expectation values
of these observables are  cohomology classes of the space
of imbeddings of a circle. One of the
resulting theories discussed in the paper has only trivalent interactions and, 
irrespective of the actual
dimension,  looks like a $3$\ndash dimensional Chern--Simons theory.
\end{abstract}
\maketitle

\section{Introduction}
Knot invariants can be obtained as expectation values of Wilson loops (i.e.,
traces of holonomies) in Chern--Simons theory \cite{W}.

The same result can be obtained in the $3$\ndash dimensional $BF$ 
theory ``with a cosmological term'' \cite{W2,CCFM}.

The nice feature of $BF$ theory \cite{S,H},
as opposed to Chern--Simons theory, is that it
can be defined in any dimension and always has a quadratic term around which
one can start a perturbative expansion.

On the other hand, apart from the $3$\ndash dimensional case 
\cite{CCM,CCFM,C1,C2}, 
it is rather difficult to find nontrivial observables for $BF$ theories
(see \cite{CM,CCFM,CCR,FMZ} for the search of surface observables in the 
$4$\ndash dimensional case).

However, the formulae for the perturbative expansion of the expectation
value of a Wilson loop in the $3$\ndash dimensional Chern--Simons theory
(e.g., in the approach of \cite{BT}) allow for a natural generalization
in any dimension \cite{CCL}. The invariants defined in this way are
cohomology classes of Vassiliev finite type \cite{V}
on the space of imbeddings of $S^1$ into $\bbR^n$. As in the 
$3$\ndash dimensional
case, they are moreover related to certain ``graph cohomologies,'' as 
originally suggested by Kontsevich \cite{K}.

This led us to look for loop observables for
$BF$ theories in any dimension whose expectation values should
yield the above invariants.

In this paper we introduce these observables   ``on shell'' (i.e., 
upon using the equation of motions) and describe some results
about their off-shell extension
(which relies on the use of the Batalin--Vilkovisky \cite{BV} formalism).
For all technical details, as well as for the complete proofs of the four 
Theorems
contained in this paper, we refer to \cite{CR}.

A very interesting feature of the simplest of these observables (as well
as of the cohomology classes described in \cite{CCL}) is that, despite of
the dimension, they are always {\em as if we were dealing with the $3$\ndash dimensional
$BF$ theory with a cosmological term}, which is tantamount to considering the
Chern--Simons theory. This is our interpretation of Witten's ideas 
about a Chern--Simons theory for strings \cite{W-1}.

\begin{Ack}
We thank R.~Longoni for carefully reading the manuscript and
for many useful discussions.
\end{Ack}

\section{$BF$ theories}
The $n$\ndash dimensional
$BF$ theory is a topological quantum field theory defined 
in terms of a connection $1$\ndash form $A$ over a principal 
$G$\ndash bundle $P\to M$,
with $\dim M=n$, and a tensorial $(n-2)$\ndash form $B$ of the 
$\ad$ type.

Following the ideas of \cite{CCR},
we begin with a geometrical description of gauge invariant
functionals of $A$ and $B$ taking values in the free loop space $\LL M$.

The basic functional is of course the Wilson loop $\Tr_\rho\hol$. Here
we follow the following {\em nonstandard}
\begin{Conv*}
$\hol(\gamma;A)$ denotes the group element associated to
the $A$-parallel transport along the loop
$\gamma:[0,1]\to M$
from $\gamma(1)$ to $\gamma(0)$. (The usual holonomy is the inverse
of our $\hol$.)
\end{Conv*}
We then associate to each smooth loop $\gamma$ in $M$ 
and to any connection $A$ the $A$\ndash horizontal lift of 
$\gamma$: $[0,1]\ni t\mapsto \gamma_A(t)$.

If we  saturate the 
$(n-2)$\ndash form $B$  with
the tangent vector $\dot \gamma_A (t)$, we obtain an $(n-3)$\ndash form
$B(\dot\gamma_A(t))$ defined along the horizontal path $\gamma_A$ and hence 
an $(n-3)$\ndash form on $LM$, depending on the connection $A$.   

Then, for any representation $\rho$ of $G$, 
the following object is well-defined and gauge invariant
in any dimension:
\begin{equation}
h_{k,\rho}(\gamma;A,B)= 
\underset{0<t_1<\ldots<t_k<1}\int
\Tr_\rho \left[B(\dot\gamma_A(t_1))\wedge \ldots\wedge 
B(\dot\gamma_A(t_k))
\hol(\gamma; A)\right].
\lab{PB}
\end{equation}
If  $n=3$, then $A+\kappa\,B$ is also a connection  and  
\eqref{PB} is equal to the $k$\ndash th Taylor coefficient in the 
$\kappa$\ndash expansion of $\Tr_{\rho}\hol(\gamma,A+\kappa\,B)$.
When $n>3$, then \eqref{PB} is a differential form over $LM$ of degree
$k(n-3)$. 

Unfortunately, however, \eqref{PB}
is not a good observable for the $BF$ theory 
since it is not invariant under the 
full set of its symmetries, see below \eqref{delta}.

If we modify \eqref{PB}
so as to get a ``good observable,'' then 
the relevant vacuum expectation values
will produce elements of the $k(n-3)$\ndash rd cohomology of knots 
(imbedded loops) as explained in the following. 
Observe that in the case $n=3$ we recover the usual Vassiliev knot invariants.



\subsection{Action functional and symmetries of $BF$ theory}
We write the action functional of the $BF$ theory under the following
\begin{Ass}
We assume that the Lie algebra
$\frg$ of $G$
possesses a nondegenerate $\Ad$\ndash invariant bilinear form 
$\braket{}{}$. 
(For example, if $\frg$ is semisimple, we may take the Killing form.)
\lab{ass-nondeg}
\end{Ass}

We extend this form to $\Omega^*(M,\ad P)$ in the usual way.
Then we define
\begin{equation}
S := \int_M \braket B{F_A},
\lab{SBF}
\end{equation}
where $F_A$ is the curvature $2$\ndash form of $A$.

The Euler--Lagrange equations of motion read
\[
F_A=0,\qquad \dd_A B=0;
\]
thus, classically, $BF$ theory in $n$ dimensions describes flat connections 
together with covariantly closed $(n-2)$\ndash forms.
Observe that ``on shell'' (i.e., on the subspace of solutions)
the covariant derivative is a coboundary operator.

\begin{Ass}
We will assume in the following that $M$ is a compact manifold, that $P$
is a trivial bundle and that
there is a flat connection $A_0$ on $P$ such that
all cohomology groups $H^*_{\dd_{A_0}}(M,\ad P)$
are trivial. 

Moreover, by ``on shell'' we will always mean ``on 
the subspace of those $(A_0,B_0)$ which satisfy the equations of motion with
$A_0$ of this kind.''
\lab{ass-triv}
\end{Ass}

The symmetries under which the action is invariant correspond to an action
of the group $\calG\rtimes_{\Ad}\Omega^{n-3}(M,\ad P)$, where
$\calG$ is the group of automorphisms of $P$ (i.e., ordinary gauge 
transformations).
Infinitesimally we have
\begin{equation}
\begin{split}
\delta A &= \dd_A \xi,\\
\delta B &= \Lie B\xi + \dd_A\chi,
\end{split}
\lab{delta}
\end{equation}
with $(\xi,\chi)\in\Omega^{0}(M,\ad P)\rtimes_{\ad}\Omega^{n-3}(M,\ad P)$.

Observe that on shell this symmetries are ``reducible;'' i.e.,
a covariantly closed $\chi$ (which under our triviality assumption
is of the form $\chi=\dd_{A_0}\sigma$) acts trivially.

More precisely, each point
$(A_0,B_0)$ in the subspace of solutions has an isotropy group 
consisting of all covariantly closed $(n-3)$\ndash forms of the $\ad$ type.

With our triviality assumption, the isotropy group at each point is then
isomorphic to the group 
$\Omega^{n-4}(M,\ad P)/\dd_{A_0}\Omega^{n-5}(M,\ad P)$.

Of course, if $n>5$, also in this quotient there are nontrivial isotropy groups
which are all isomorphic to 
$\Omega^{n-5}(M,\ad P)/\dd_{A_0}\Omega^{n-6}(M,\ad P)$, and so on until we 
arrive at $\Omega^0(M,\ad P)$ which acts freely on  $\Omega^1(M,\ad P)$.

In order to consistently gauge-fix all the symmetries, one has then to resort
to the extended BRST formalism (i.e., introduce a hierarchy of ghosts for
ghosts). 

An additional problem is due to the fact that the isotropy groups
are different off shell. So, in order to work in the Lagrangian formalism which
is better suited for the perturbative expansions, one has to rely on the
whole Batalin--Vilkovisky (BV) machinery as explained in Section~\ref{sec-BV}. 

It is known that the partition function of $BF$ theory
is   related to the analytic torsion
of $M$ (see \cite{S} for the abelian case and \cite{BlT} for the
non-abelian one).

In order to get other topological invariants, one has to find
interesting BV closed observables. In the rest of this paper we will discuss
some of them, leaving most of the technical details to \cite{CR}.

Before starting with the general discussion we recall
the $3$\ndash dimensional
case as studied in \cite{CCM,CCFM}.

\subsection{The $3$\ndash dimensional $BF$ theory}
Since $B$ is a $1$\ndash form now, one can add to the pure $BF$ action
the so-called cosmological term,
\[
S_3 := \frac16\int_M \braket B{\Lie BB},
\]
and define
\[
S_\kappa:= S + \kappa^2\,S_3,\qquad\kappa\in\bbR.
\]

This action is actually equal to the difference of two Chern--Simons actions
evaluated at the connections $A+\kappa\, B$ and $A-\kappa\, B$. From this
observation one immediately gets the infinitesimal symmetries for this
theory as
\begin{equation}
\begin{split}
\delta_\kappa A &= \dd_A \xi + \kappa^2\,\Lie B\chi,\\
\delta_\kappa B &= \Lie B\xi + \dd_A\chi.
\end{split}
\lab{deltakappa}
\end{equation}

Another consequence is that, for any loop $\gamma$ and any representation
$\rho$ of the Lie algebra $\frg$ of $G$, the Wilson loops
$\Tr_\rho\hol(\gamma;A\pm\kappa\,B)$ are observables.
We then have
\begin{multline}
\calH_\rho(\kappa;\gamma;A,B) = \Tr_\rho\hol(\gamma;A+\kappa\,B)
=\Tr_\rho\hol(\gamma;A) +
\sum_{k=1}^\infty \kappa^{k}\,h_{k,\rho}(\gamma;A,B),
\lab{hk}
\end{multline}
where $h_{k,\rho}(\gamma;A,B)$ is the $k$\ndash th Taylor coefficient
in the $\kappa$\ndash expansion of $\Tr_\rho\hol(\gamma;A+\kappa\,B)$
introduced in \eqref{PB}.

By the previous discussion it follows clearly that both the even and the odd
part of $\calH$ are observables.

In order to give an explicit description of $h_{k,\rho}$, it is better to view
$\gamma$ as a periodic mapping from $[0,1]$ to $M$. We denote by
$H(\gamma;A)|_{s}^{t}$ the group element determined, in a given trivialization,
by the $A$\ndash parallel transport along $\gamma$
from the point $\gamma(t)$ to the point $\gamma(s)$. 
Then we can rewrite \eqref{PB} in the form
\begin{multline*}
h_{k,\rho}(\gamma;A,B)= \int_{\triangle_k}\Tr_\rho\Big[ H(\gamma;A)|_0^{t_1}
\,B_1\,H(\gamma;A)|_{t_1}^{t_2}\,
B_2\cdots\\
\cdots B_{k-1}\ H(\gamma;A)|_{t_{k-1}}^{t_k}\,B_k\,
H(\gamma;A)|_{t_k}^{1}\Big],
\end{multline*}
where $\triangle_k$ is the $k$-simplex $\{0<t_1<t_2<\dots<t_k<1\}$
and $B_i$ is a shorthand notation for the pullback of $B$ via the map
$\gamma_i:\triangle_k\to M$, $\gamma_i(t_1,\dots,t_k)=\gamma(t_i)$.

Observe that each $h_{k,\rho}(\bullet;A,B)$ defines a function on the loop 
space $\LL M$ of $M$. It is not difficult to check that---modulo the
equations of motion
\[
F_A+\frac{\kappa^2}2\,\Lie BB=0,
\qquad
\dd_A B=0,
\]
of $S_\kappa$---the functions $\calH(\bullet;A,B)$ are locally constant
on $\LL M$. 

Quantization then requires a regularization, viz., point splitting or, in a more
precise formulation, the blowing up of the diagonals of configuration spaces.
So one has to consider imbeddings,
instead of generic loops, and to introduce a framing.
The expectation values of the $\calH$s
then define
locally closed functions on the space $\Imb_f(S^1, M)$ of framed imbeddings,
i.e., framed knot invariants.


\subsection{A first glimpse to higher dimensions}\lab{ssec-ghd}
A straightforward generalization of \eqref{hk} to higher dimensions 
exists as discussed at the beginning of this section.
It yields
forms on $\LL M$ instead of functions since only one form degree for each
field $B$ is saturated by the integration. So now
\[
h_{k,\rho}(\bullet;A,B)\in \Omega^{(n-3)k}(\LL M).
\]

Assume now that $A$ is a flat connection and that $B$ is covariantly closed, as
classical solutions of the pure $BF$ theory are.
Then it is not difficult to check that, for any odd $n>3$, $h_{k,\rho}$ is 
closed.
This follows from the generalized Stokes theorem.

More precisely, $h_{k,\rho}$ is a form on $\Hat M=\LL M\times\calA\times
\Omega^{n-2}(M,\ad P)$, where $\calA$ is the space of connections on $M$,
and we restrict it to the subspace where $F_A=\dd_A B=0$.

Let $\Hat\dd $ be the differential on $\Hat M$. 
In the computation of $\Hat\dd h_{k,\rho}$ we can switch the integral
with the differential. We then get $\Hat\dd$ acting on a function 
$\eta_{k,\rho}$ on $\Hat M\times\triangle_k$:
\[
\Hat\dd h_{k,\rho}=\int_{\triangle_k}\Hat\dd\eta_{k,\rho}.
\] 

Let $\dd_k$ be the differential on $\triangle_k$ and $\dd:=\Hat\dd\pm\dd_k$
the differential on $\Hat M\times\triangle_k$. By adding and subtracting
$\dd_k$ we get then
\[
\Hat\dd h_{k,\rho} = \int_{\triangle_k}\dd \eta_{k,\rho}
\mp \int_{\triangle_k}\dd_k \eta_{k,\rho}.
\]
The first term on the r.h.s.\ is then easily seen to vanish in our
hypotheses $F_A=\dd_A B=0$ since
\[
\dd H(\gamma;A)|_{s}^{t} = -A(s)\, H(\gamma;A)|_{s}^{t}
+  H(\gamma;A)|_{s}^{t}\,A(t)
\]
when $A$ is flat.

The second term can then be computed using the Stokes theorem. The
codimension-one boundaries of the simplex corresponding to
the collapse of consecutive points yield terms containing $B^2$ which
vanishes for dimensional reasons ($n>4$). The remaining codimension-one
boundaries correspond to $t_0=0$ or $t_k=1$. It is not difficult to check, 
using the cyclic property of the trace, that these terms cancel each other.

If the dimension $n$ is even, $n>4$, the only problem in the above discussion
arises at the last step since, for $k$ even, the two terms coming from
$t_0=0$ and $t_k=1$ sum up instead of canceling each other.
For $k$ odd however everything works as before.

Similar computations allow to show that the $h_{k,\rho}$ are invariant (modulo
exact forms on $\LL M$) under the
symmetries \eqref{delta} either if $n$ is odd and greater then 5 or
if $n$ is even and greater than 4 and $k$ is odd.

In sections~\ref{sec-odd} and~\ref{sec-gen} we will describe how this 
discussion can be extended ``off 
shell'' and how the cases $n=4$ and $n=5$ will be included.


\section{The BV quantization of $BF$ theories}\lab{sec-BV}
\subsection{The BRST operator}
In order to deal with the symmetries \eqref{delta} of \eqref{SBF}
in the functional integral, one has to introduce the BRST operator
\begin{equation}
\begin{split}
\dBRST A &= \dd_A c,\\
\dBRST B &= \Lie B c + \dd_A\tau_1,
\end{split}
\lab{deltaBRST}
\end{equation}
where $c$ and $\tau_1$ are ghosts, i.e., forms on the space of fields with 
values
in $\Omega^0(M,\ad P)$ and $\Omega^{n-3}(M,\ad P)$ respectively.
As usual in gauge theories one also defines
\[
\dBRST c = -\frac12\Lie cc.
\]

Because of the on-shell reducibility, one has then to introduce ghosts for
ghosts $\tau_k$ with values in $\Omega^{n-2-k}(M,\ad P)$, $k=1,\dots,n-2$, with
ghost number equal to $k\mod2$ and extended BRST operator
\[
\begin{split}
\dBRST\tau_k &= (-1)^k\, \Lie{\tau_k}c + \dd_A\tau_{k+1},\qquad k=1,\dots,n-3,\\
\dBRST\tau_{n-2} &= (-1)^n \Lie{\tau_{n-2}}c.
\end{split}
\]

\subsection{The BV formalism}
It is not difficult to check that $\delta_{\mathrm{BRST}}^2=0\mod F_A$.

The Batalin--Vilkovisky method allows then for the construction of
a nilpotent operator $\dBV$ that extends $\dBRST$ off shell.

To do so, one first introduces a partner $\phi_\alpha^+$
with values in $\Omega^*(M,\ad P)$ 
for any field or ghost $\phi^\alpha=A,B,c,\tau_1,\dots,\tau_{n-2}$ with the 
following rules:
\begin{itemize}
\item The ghost number of $\phi_\alpha^+$ is minus the ghost number of
$\phi^{\alpha}$, minus one.
\item The form degree of $\phi_\alpha^+$ is $n$ minus the form degree
of $\phi^\alpha$.
\end{itemize}
\begin{Rem}
In general, the antifields are dual to the corresponding fields.
Of course some isomorphisms may be used to identify certain spaces.
For example here we have preferred to identify the Lie algebra $\frg$
with its dual (using the bilinear from $\braket{}{}$)
so that also the antifields take values
in the space of tensorial forms of the $\ad$ type.
This will be particularly useful, e.g., in equations \eqref{defsfB} and
\eqref{defsfA}.

In the original formulation of Batalin and Vilkovisky \cite{BV} one also
identifies forms of complementary degree using a Hodge operator.
Since we do not want to introduce a metric here, we prefer to avoid this 
identification. As a consequence,  our BV antibracket \eqref{BVbracket}
will be of the form described in \cite{CF} instead of the original one.
\end{Rem}


\begin{Rem}[Sign convention]
We follow here the usual convention for the sign rules related
to the double grading given by the form degree $\deg$ and the ghost
number $\gh$.

Namely, in the case of homogenous forms $\alpha$ and $\beta$
of the $\ad$\ndash type we have
\[
\Lie{\alpha}{\beta}=-(-1)^{\deg\alpha \deg\beta+ \gh\alpha \gh\beta}\Lie{\beta}{\alpha}.
\]
Moreover, in the case of homogenous forms $\alpha$ and $\beta$ taking values
in a commutative algebra (e.g., $\bbR$) we have
\[
\alpha\wedge \beta=(-1)^{\deg\alpha \deg\beta+ \gh\alpha \gh\beta}\beta \wedge \alpha.
\]
\end{Rem}

Next we define the BV bracket of two functionals $F$ and $G$.
We use throughout Einstein's convention over repeated indices and set
\begin{equation}
\BV FG := 
\int_M\braket{F\dleft{\phi^\alpha}}{\dright{\phi_\alpha^+}G}
-(-1)^{\deg\phi^\alpha\,(n+1)}\,
\braket{F\dleft{\phi_\alpha^+}}{\dright{\phi^\alpha}G},
\lab{BVbracket}
\end{equation}
where the left and right functional derivatives are given by the following 
formula:
\[
{\dod t}{\Big|_{t=0}} F(\phi^\alpha+t\,\rho^\alpha)=
\int_M\braket{\rho^\alpha}{\dright{\phi^\alpha}F}=
\int_M \braket{F\dleft{\phi^\alpha}}{\rho^\alpha};
\]
we proceed similarly for the antifields.

As usual, the space of functionals with BV bracket is a Gerstenhaber 
algebra \cite{G}.

Finally the BV operator is defined by
\begin{equation}
\dBV := \BV\SBV{\ }
\lab{defdBV}
\end{equation}
where $\SBV$ is a solution of the master equation
\[
\BV\SBV\SBV=0
\]
such that $\SBV_{|_{\phi_\bullet^+=0}}=S$.

\subsection{The BV action for $BF$ theories}
The BV action $\SBV$ corresponding to the $BF$ action \eqref{SBF} 
can be written as
\begin{equation}
\SBV = \int_M \dbraket{\sfB}{\sfF_\sfA}
\lab{SBV}
\end{equation}
where the notations are as follows:
\begin{itemize}
\item The {\em dot product}\/ is just the wedge product between forms taking
values in an associative algebra---e.g., $\bbR$ or $\frg$ itself if it
associative as in Section~\ref{sec-gen}---but with a shifted 
degree;
viz., for two homogenous forms $\alpha$ and $\beta$ we set
\[
\alpha\cdot\beta:=(-1)^{\gh\alpha\,\deg\beta}\,\alpha\wedge\beta,
\]
where $\gh$ denotes the ghost number.

We extend then the bilinear form $\braket{}{}$ to forms with shifted degree
by setting
\[
\dbraket\alpha\beta := (-1)^{\gh\alpha\,\deg\beta}\,
\braket\alpha\beta.
\]
Similarly we define the {\em dot Lie bracket}\/ for two homogeneous forms of the 
$\ad$ type by
\[
\dLie\alpha\beta:=(-1)^{\gh\alpha\,\deg\beta}\,\Lie\alpha\beta.
\]
These definitions are then extended by linearity. 

An easy check shows
that the dot product in the case of a commutative algebra and
the dot Lie bracket are, respectively, a graded commutative product and
a graded Lie bracket with respect to a new grading called
the {\em total degree}\/ that is defined as the form degree plus the ghost 
number. Moreover, $\dd_{A_0}$ is still a differential for 
the dot algebras.

\item The ``super $B$\ndash field'' $\sfB$ is defined by
\begin{equation}
\sfB = \sum_{k=1}^{n-2} (-1)^{\frac{k(k-1)}{2}}\tau_{k}+B+(-1)^n A^++c^+
\in\Omega^*(M,\ad P),
\lab{defsfB}
\end{equation}
and has total degree equal to $n-2$.
\item The ``supercurvature'' $\sfF_\sfA$ of the ``superconnection''
\begin{equation}
\sfA = (-1)^{n+1} c+A+(-1)^n B^++\sum_{k=1}^{n-2} (-1)^{n(k+1)+\frac{k(k-1)}{2}}\tau_{k}^+
\lab{defsfA}
\end{equation}
is given by the usual formula. In order to write it down, it is better
to choose a background connection $A_0$ and to define the tensorial form
$\sfa=\sfA-A_0$ of total degree one. Then
\[
\sfF_\sfA = F_{A_0} + \dd_{A_0}\sfa + \frac12\, \dLie\sfa\sfa.
\]
In general we will choose $A_0$ to be a flat connection
as in Assumption~\ref{ass-triv}.
\item By $\int_M$ we then mean the integral of all the terms of form degree
equal to $n$. Observe that as a consequence $\SBV$ has then ghost number zero.
\end{itemize}

\begin{Rem}
We may observe that there is a superspace formulation of
\eqref{SBV} obtained
by introducing superpartners to the coordinates of $M$ and redefining
$\sfA$ and $\sfB$ accordingly. In this way we would follow the pattern
described in \cite{DG}.

Special cases were already discussed in \cite{C2,CF} (two dimensions) and
\cite{CRR,BF7} (four dimensions). See also \cite{W3,AKSZ} for the case of the
$3$\ndash dimensional Chern--Simons theory. 

For later purposes---viz., in order to define loop observables as in the
following sections---it is however better to work in our setting.
\end{Rem}

We have the following general result \cite{CR}:
\begin{Thm}
The action $\SBV$ satisfies the master equation in any dimension.
\lab{thm-me}
\end{Thm}

We conclude this section by giving the explicit action of the BV operator
\eqref{defdBV} on the ``superfields'' $\sfA$ and $\sfB$. In order to give
neater formulae, it is better to define a new BV operator with shifted
degree; viz., for a homogeneous form $\alpha$, we set
\[
\sfdelta\alpha := (-1)^{\deg\alpha}\,\dBV\alpha.
\]
One can show that $\sfdelta$ is a differential for the dot algebras and
that it anticommutes with $\dd_{A_0}$.

Then we obtain \cite{CR}
\begin{equation}
\begin{split}
\sfdelta\sfA &= (-1)^n\,\sfF_\sfA,\\
\sfdelta\sfB &= (-1)^n\,\dd_\sfA\sfB,
\end{split}
\lab{sfdelta}
\end{equation}
with
\[
\dd_\sfA\sfB = \dd_{A_0}\sfB + \dLie{\sfa}\sfB.
\]

Upon using the above equations, we can then prove Thm.~\ref{thm-me}
by simply checking that $\sfdelta S_{\mathrm{BV}}=0$,  
as follows from the the $\ad$\ndash invariance of $\braket{}{}$
and from the Stokes theorem.

\subsection{The BV Laplace operator and the BV observables}
In the quantum version of the BV formalism---i.e., when dealing with 
functional integrals with weight $\exp(\ii/\hslash)S$---one has then to 
introduce the so-called BV Laplace operator
$\DBV$ and to verify that the {\em quantum}\/ master equation
\[
\BV\SBV\SBV - 2\ii\hslash \DBV\SBV=0
\]
is satisfied.

The very definition of $\DBV$ 
relies on a regularization of the theory, which we do not discuss here.
We only recall the formal properties of
$\DBV$; viz.:
\begin{enumerate}
\item $\DBV$ is a coboundary operator on the space of functionals;
\item for any two functionals $F$ and $G$,
\begin{equation}
\DBV(F\,G) = (\DBV F)\,G + (-1)^{\gh F}\,F\DBV G+
(-1)^{\gh F}\BV FG.\lab{DBVFG}
\end{equation}
\end{enumerate}
The space of functionals with the BV bracket and the BV Laplacian is a
so-called BV algebra.

To give an explicit definition of the BV Laplacian one has to introduce
some extra structures (e.g., a Riemannian metric on $M$) and a regularization.
The main property however is that the BV Laplace operator contracts
each field with the Hodge dual of the
corresponding antifield at the same point in $M$.

Under our assumptions,
one can then prove \cite{CR} that $\DBV\SBV=0$ for $\SBV$ in \eqref{SBV}. 

So $\SBV$ is also a solution
of the quantum master equation. This implies that its partition function
is independent of the choice of gauge fixing.

A consequence of the properties of the BV operators is that the operator
\[
\OBV := \dBV -\ii\hslash\DBV
\]
is a coboundary operator if{f} $\SBV$ satisfies the quantum master equation.

The main statement in the BV formalism is that the $\OBV$\ndash cohomology
of ghost number zero yields all the meaningful observables. 
More precisely, this means that
the expectation value of an $\OBV$\ndash closed functional is independent of the
gauge fixing and that the expectation value of an $\OBV$\ndash exact functional
(or of a functional of ghost number different from zero) vanishes.

In the next sections we will discuss some BV observables of $BF$ theories
associated to $\LL M$ (or better to $\ImbfM$). To do so, it is however better
to use shifted versions of the operators $\OBV$ and $\DBV$ as well, viz.:
\begin{align*}
\sfOmega\alpha&:=(-1)^{\deg \alpha}\OBV \alpha ,\\
\sfDelta\alpha&:=(-1)^{\deg \alpha}\DBV \alpha.
\end{align*}

\section{Generalized Wilson loops in odd dimensions}\lab{sec-odd}
At this point we are ready to define the correct generalization of
the observables $\calH$ defined in \eqref{hk}.

Formally the new observable is still the trace of the ``holonomy
of $\sfA+\kappa\,\sfB$''
\begin{equation}
\sfcalH_\rho(\kappa;\sfA,\sfB):=
\Tr_\rho \sfhol(\sfA+\kappa\,\sfB)\in\Omega^*(\LL M),
\lab{defcalH}
\end{equation}
where the ``holonomy'' $\sfhol$ is now defined in terms of iterated integrals
as follows:
First we write $\sfA=A_0+\sfa$. Then we set
\[
\Tr_\rho \sfhol(\sfA+\kappa\,B):=
\Tr_\rho \hol(A_0)+
\sum_{l=1}^\infty h_{l,\rho}(A_0,\sfa+\kappa\,\sfB),
\]
where
\begin{multline*}
h_{l,\rho}(A_0,\sfa+\kappa\,\sfB) = 
\int_{\triangle_{l}}\Tr_\rho\Big[
H(A_0)|_0^{t_1}\cdot (\sfa_1+\kappa\,\sfB_1)
\cdot H(A_0)|_{t_1}^{t_2}\cdot\\
\cdot(\sfa_2+\kappa\,\sfB_2)
\cdots H(A_0)|_{t_{l-1}}^{t_l}\cdot (\sfa_l+\kappa\,\sfB_l)\cdot
H(A_0)|_{t_l}^{1}\Big].
\end{multline*}
Here $\sfa_i$ and $\sfB_i$ are shorthand notations
for the pullbacks of $\sfa$ and $\sfB$ via 
$\mathrm{ev}_i:\LL M\times\triangle_{l}
\to M$, $(\gamma;t_1,\dots,t_l)\mapsto\gamma(t_i)$. Moreover, 
the $H(A_0)|_{s}^{t}$'s denote the group elements associated to
parallel transports as functions on $\LL M$.

\begin{Rem}
The above integrals should be better viewed as integrations along the
fiber of the trivial bundles $\triangle_l\times\LL M\to \LL M$.
That is, the integrals are zero whenever the form degree is less
than the dimension of the simplex and yield a form on $\LL M$ whenever
the form degree exceeds the dimension of the simplex.

Also observe that ${\sfcalH}$ is a sum of terms with
different ghost number and different form degree on $\LL M$.
\end{Rem}

Of course, we cannot expect $\sfcalH$ to be an observable for the pure $BF$
theory. We can however consider a fake ``higher dimensional $BF$ theory
with cosmological term'' as follows: We first define
\begin{equation}
\sfS_3(\sfB):= \frac16\int_M \dbraket\sfB{\dLie \sfB\sfB},
\lab{S3}
\end{equation}
where again we consider only the terms of form degree equal
to $n$, so $\sfS_3$ has ghost number $2(n-3)$.
Then we consider the functional
\begin{equation}
\widehat{\sfcalH}_\rho(\hslash,\kappa;\sfA,\sfB):=
\left\{ \exp[(\ii/\hslash)\,\kappa^2\,\sfS_3(\sfB)]\cdot
\sfcalH_\rho(\kappa;\sfA,\sfB)\right\}_0,
\lab{hHoddn}
\end{equation}
where $\{\}_0$ means taking the terms with ghost number zero.

The functional $\widehat{\sfcalH}_\rho$ is well-defined for
any loop in $M$. However, in order to avoid problems with
the BV Laplace operator $\DBV$, we must restrict ourselves to the space
of framed imbeddings $\ImbfM$.

\begin{Thm}
For any $\kappa$ and $\rho$ and any odd dimension $n$,
$\widehat{\sfcalH}$, as a functional taking values in the
forms on $\ImbfM$,
is $\OBV$\ndash closed modulo $\dd$\ndash exact forms and $\dd$\ndash closed 
modulo $\OBV$\ndash exact terms.
In other words, $[\widehat{\sfcalH}]$ is an $H^*(\ImbfM)$\ndash valued
observable.
\lab{thm-Hoddn}
\end{Thm}

\begin{proof}[Proof (Sketch)]
The main idea of the proof relies on the identity
\begin{multline}
\sfOmega\left\{ \exp[(\ii/\hslash)\,\kappa^2\,\sfS_3(\sfB)]\cdot
\sfcalH_\rho(\kappa;\sfA,\sfB)\right\}
=\exp[(\ii/\hslash)\,\kappa^2\,\sfS_3(\sfB)]\cdot
\boldsymbol{\delta}_{\kappa} \sfcalH_\rho(\kappa;\sfA,\sfB),
\lab{sfdeltakappa}
\end{multline}
where $\boldsymbol{\delta}_{\kappa}$ is the following coboundary operator:
\[
\begin{split}
\boldsymbol{\delta}_\kappa\sfA &= -\sfF_\sfA -
\frac{\kappa^2}2\,\dLie\sfB\sfB,\\
\boldsymbol{\delta}_\kappa\sfB &= -\dd_\sfA\sfB.
\end{split}
\]
Observe that for $n\not=3$, $\boldsymbol{\delta}_{\kappa}$ is a differential only
for the $\bbZ_2$\ndash reduction of the graded algebra of functionals.

Using $\boldsymbol{\delta}_{\kappa}$ is like working with a cosmological term,
and, upon using the generalized Stokes theorem, one gets
\[
(\dd+\boldsymbol{\delta}_{\kappa})\,\sfcalH_\rho(\kappa;\sfA,\sfB)=0,
\]
which proves the theorem.

Equation \eqref{sfdeltakappa} is a consequence of \eqref{DBVFG}
and of the following identities:
\[
\sfdelta\sfS_3=0, \quad
\sfDelta\exp[(\ii/\hslash)\,\kappa^2\,\sfS_3(\sfB)]=0, \quad  
\sfDelta\sfcalH_\rho =0.
\]
The first identity follows from \eqref{sfdelta}, from the fact that
$\braket{}{}$ is $\ad$\ndash invariant and from Stokes theorem.

The second identity holds since $\sfS_3$ depends only on $\sfB$ and as a 
consequence of the already discussed property according to which
the BV Laplace operator contracts
each field with the Hodge dual (for some Riemannian metric on $M$) of the
corresponding antifield at the same point in $M$.

The last identity is ``formally'' (i.e., modulo regularization problems
for $\DBV$) true if $\gamma$ does not have transversal self-intersections,
for the same reason as above.
However, in order to rely upon this last identity confidently, we must then
restrict ourselves to framed imbeddings and put each component of $\sfA$
on the imbedding and each component of $\sfB$ on its companion (as done
in \cite{CCM}).
\end{proof}

As a consequence, the expectation value of $\widehat{\sfcalH}$ 
is (up to anomalies) a cohomology class on the space of (framed)
imbeddings of $S^1$ into $M$.



\begin{Rem}
If we set all the antifields to zero, $\widehat{\sfcalH}$ reduces to a 
sum of $h_{k,\rho}(A,B)$'s; so it is the off-shell generalization
we were looking for in subsection~\ref{ssec-ghd}.

Moreover,  the expectation value of $\widehat{\sfcalH}$ w.r.t.\ the
pure $BF$ theory is in three dimensions the same as the expectation
value of $\calH$ in the $BF$ theory with cosmological term.
\end{Rem}

\begin{Rem}
$\widehat{\sfcalH}$ is a genuine quantum observable
since its limit for $\hslash\to0$ is not defined.

However, one might replace $\kappa$ with $\hslash\kappa$. In this way,
$\widehat{\sfcalH}$ becomes a formal power series in $\hslash$. 
The zeroth order
term is just $\Tr_\rho\hol(A)$.

Since $\widehat{\sfcalH}$ is an observable for any $\kappa$, so is
its odd part (in $\kappa$) $\widehat{\sfcalH}^o$. It is not difficult to see
that $\widehat{\sfcalH}^o_\rho(\hslash,\hslash\kappa;\sfA,\sfB)/\hslash$
is as well a formal power series in $\hslash$ and that its zeroth-order term is
the observable 
\[
\sfh_{1,\rho}(\sfA,\sfB):=
\dod\kappa \Big|_{\kappa=0}
h_{1,\rho}(A_0,\sfa+\kappa B),
\]
which is the off-shell extension of $h_{1,\rho}(A,B)$.

Therefore, the observables
$\widehat{\sfcalH}_\rho(\hslash,\hslash\kappa;\sfA,\sfB)$ 
and $\widehat{\sfcalH}^o_\rho(\hslash,\hslash\kappa;\sfA,\sfB)/\hslash$ are
nontrivial quantum deformations of, respectively, the ordinary Wilson loop and 
$\sfh_{1,\rho}$.
\end{Rem}

\section{Other loop observables}\lab{sec-gen}
In order to generalize some of the results of the previous section
and in order to define more general observables we make the following
\begin{Ass}
We assume that the Lie algebra $\frg$ is obtained from an associative
algebra with trace $\Tr$
(e.g., we may take $\frg=\frgl(N)$ with the usual trace of matrices). 
In this case, we assume that our $\ad$\ndash invariant
bilinear form $\braket\xi\eta$ is given by $\Tr(\xi\,\eta)$ and,
according to Assumption~\ref{ass-nondeg}, we further assume that 
it is nondegenerate. 
Moreover, we consider only
representations $\rho$ of $\frg$ as an associative algebra.
\lab{ass-ass}
\end{Ass}

\subsection{Generalized Wilson loops in even dimensions}
The observable defined in the previous section does not work in even
dimensions essentially because the ``cosmological term'' $\sfS_3$
vanishes when $\sfB$ has even total degree. 

We can cure this problem thanks to Assumption~\ref{ass-ass} by defining
instead
\[
\sfO_3 (\sfB):=\frac13\int_M \Tr(\sfB\cdot\sfB\cdot\sfB).
\]

We have already seen in subsection~\ref{ssec-ghd}
that the even part of $\calH$ does not work. So we consider only
\begin{equation}
\widehat{\sfcalH}^o_\rho(\hslash,\kappa;\sfA,\sfB):=
\left\{ \exp[(\ii/\hslash)\,\kappa^2\,\sfO_3(\sfB)]\cdot
\sfcalH^o_\rho(\kappa;\sfA,\sfB)\right\}_0,
\lab{hHevenn}
\end{equation}
where $\sfcalH^o$ is the odd part of $\sfcalH$, which is defined exactly
as in  the odd-dimensional case---see equations \eqref{defcalH} and following.

We have then the following analogue (with analogous proof) 
of Thm.~\ref{thm-Hoddn}:
\begin{Thm}
For any $\kappa$ and $\rho$ and any even dimension $n$,
$[\widehat{\sfcalH}^o]$ is an $H^*(\ImbfM)$\ndash valued observable.
\lab{thm-Hevenn}
\end{Thm}

\begin{Rem}
Similarly to what happens in the odd-dimensional case, 
$\widehat{\sfcalH}^o$ reduces
to a sum of $h_{2k+1,\rho}$ as the antifields are set to zero.

Moreover,
$\widehat{\sfcalH}^o_\rho(\hslash,\hslash\kappa;\sfA,\sfB)/\hslash$
is still a nontrivial quantum deformation of $\sfh_{1,\rho}$.

However, we do not find in even dimensions a  nontrivial quantum deformation of
the ordinary Wilson loop $\Tr_\rho\hol$.
\end{Rem}

\subsection{Loop observables with more then cubic interactions}
The ``cosmological terms'' $\sfS_3$ and $\sfO_3$ give rise, in the perturbative
expansion, to trivalent vertices.

If we work with Assumption~\ref{ass-ass}, we can define more general 
interaction terms:
\[
\sfO_r(\sfB) := \frac1r\int_M \Tr\sfB^r.
\]
Observe that in odd dimensions $\sfO_{r}$ vanishes if $r$ is even.

Next, for any two given sequences $\bmu=\{\mu_1,\mu_2,\dots\}$ 
and $\blambda=\{\lambda_1,\lambda_2,\dots\}$, we define
\begin{multline}
\widetilde{\sfcalH}_\rho(\hslash,\bmu,\blambda;\sfA,\sfB):=
\Bigg\{ \exp\left[(\ii/\hslash)\sum_{r=1}^\infty \mu_r\,\sfO_{r+1}
(\sfB)\right]
\cdot\Tr_\rho\sfhol\left(\sfA+\sum_{s=1}^\infty\lambda_s\,\sfB^s\right)
\Bigg\}_0.
\lab{tHgeneral}
\end{multline}
We then denote by $\widetilde{\sfcalH}_\rho^o$ the odd part of
$\widetilde{\sfcalH}_\rho$ under $\blambda\to-\blambda$.

We have then
\begin{Thm}
In odd dimensions, $[\widetilde{\sfcalH}_\rho]$ is an $H^*(\ImbfM)$\ndash valued
observable whenever the following conditions are satisfied
\begin{gather*}
\mu_{2l-1}=\lambda_{2l}=0,\quad\forall l,\\
\mu_{2l}=\sum_{\substack{i,j\ge0\\i+j=l-1}}\lambda_{2i+1}\,\lambda_{2j+1},
\quad\forall l.
\end{gather*}

In even dimensions, $[\widetilde{\sfcalH}_\rho^o]$ is an $H^*(\ImbfM)$\ndash valued
observable whenever the following conditions are satisfied
\[
\mu_l=\sum_{\substack{i,j\ge1\\i+j=l}} \lambda_i\,\lambda_j,\quad\forall l.
\]
\lab{thm-tHgeneral}
\end{Thm}
The proof is a direct generalization of the proof of Thm.~\ref{thm-Hoddn}.

\section{Conclusions}
In this paper we have defined some $H^*(\ImbfM)$\ndash valued observables
for $BF$ theories on a trivial principal $G$-bundle $P\to M$ associated
to any representation of the Lie algebra $\frg$. 

Our ideas extend naturally to nontrivial bundles as well, though 
we did not consider this extension here for the sake of simplicity.

The expectation values of these observables define then classes in the
cohomology of the space of framed imbeddings of $S^1$ into $M$, which we have
assumed to be compact in order to simplify the discussion.

Of course a very interesting case is $M=\bbR^n$, which is not compact.
The only extra technical point here 
is that one has to specify the correct behavior
at infinity of all the fields and antifields. This done,
the trivial connection satisfies the hypotheses of Assumption~\ref{ass-triv}.

The perturbative expansion of the expectation values
in the case  $M=\bbR^n$ around the
trivial connection is then obtained in terms of the configuration space
integrals discussed in \cite{CCL}, where however only the framing-independent
coholomogy classes were considered explicitly.

Notice that in this paper we have {\em not}\/ defined observables with trivalent
interactions and an even number of $B$\ndash fields placed on the imbedding in the
even-dimensional case. Thus, we cannot obtain the nontrivial class of 
imbeddings represented by the diagram cocycle of Figure~4 in \cite{CCL}.
This suggests that there might exist other observables than those we have 
considered here.

On the other hand, one may use the combinatorics of this quantum field
theory to obtain new nontrivial diagram cocycles.

\end{document}